\documentclass{amsart}
\usepackage[english]{babel}
\usepackage{amsmath, amsthm, amssymb}
\usepackage{enumerate}
\usepackage{hyperref}
\usepackage{cleveref}

\usepackage[a4paper, margin=1in]{geometry}

\usepackage{tabularx}
\usepackage{caption}
\usepackage{subcaption}
\usepackage{color}
\usepackage[thinlines]{easytable}
\usepackage{physics}

\usepackage{tikz}
\usetikzlibrary{arrows,shapes}
\usetikzlibrary{arrows.meta}
\usetikzlibrary{decorations.markings}
\usetikzlibrary{intersections}
\usetikzlibrary{positioning}

\newcommand{\arrowRight}{\tikz \draw[-Stealth] (-1pt,0) -- (1pt,0);}
\newcommand{\arrowLeft}{\tikz \draw[-Stealth] (1pt,0) -- (-1pt,0);}
\newcommand{\arrowUp}{\tikz \draw[-Stealth] (0,-1pt) -- (0,1pt);}
\newcommand{\arrowDown}{\tikz \draw[-Stealth] (0,1pt) -- (0,-1pt);}

\newcommand{\arrowLeftRed}{\tikz \draw[-{Stealth[red]}] (1pt,0)--(-1pt,0);}
\newcommand{\arrowUpRed}{\tikz \draw[-{Stealth[red]}] (0,-1pt) -- (0,1pt);}
\newcommand{\arrowDownRed}{\tikz \draw[-{Stealth[red]}] (0,1pt) -- (0,-1pt);}

\title{The Six-Vertex Model with a Non-Standard Boundary Condition}
\date\today

\address[M. Dole\v{z}\'{a}lek, M. Ra\v{s}ka, E. Sgallov\'{a}, E. N. Stucky, M. Zindulka]{Charles University, Faculty of Mathematics and Physics, Department of Algebra,
Sokolovsk\'{a} 83, 186 75 Praha 8, Czech Republic}
\address[E. Sgallov\'{a}]{Michigan State University, East Lansing, USA}
\author[M. Dole\v{z}\'{a}lek]{Mat\v{e}j Dole\v{z}\'{a}lek}
\email{matej@gimli.ms.mff.cuni.cz}
\author[M. Ra\v{s}ka]{Martin Ra\v{s}ka}
\email{raska.martin@gmail.com}
\author[E. Sgallov\'{a}]{Ester Sgallov\'{a}}
\email{sgallova@msu.edu}
\author[E. N. Stucky]{Eric Nathan Stucky}
\email{stuck127@umn.edu}
\author[M. Zindulka]{Mikul\'{a}\v{s} Zindulka}
\email{mikulas.zindulka@matfyz.cuni.cz}

\thanks{We acknowledge support by Czech Science Foundation (GA\v{C}R) grant 21-00420M and the Charles University Research Centre program UNCE/SCI/022.}

\theoremstyle{plain}
\newtheorem{theorem}{Theorem}
\newtheorem{lemma}[theorem]{Lemma}

\newtheorem*{example}{Example}
\theoremstyle{definition}
\newtheorem*{definition}{Definition}

\begin{document}

\begin{abstract}
We consider the enumeration of states in the Brubaker--Bump--Friedberg six-vertex model, whose boundary conditions are determined by an integer partition. In general, we find the number of states is a polynomial in the largest part of the partition. By explicating this technique, we also enumerate the states completely for hook shapes and staircases.
\end{abstract}

\maketitle

\section{Introduction}

One celebrated generalization of permutations are the \emph{alternating sign matrices} (ASMs), which are matrices of $1$s, $0$s, and $-1$s such that in each row and column of the matrix, the sum of all entries is $1$ and the nonzero entries alternate in sign. ASMs were first noticed by Robbins and Rumsey in an analysis of Dodgson condensation (see~\cite{Bre99}). They were quickly discovered to be in bijective correspondence with many other combinatorial objects~\cite{MRR82, MRR83}. 

The problem of enumerating ASMs attracted considerable attention because despite resisting many early efforts at proof, the formula itself is remarkably simple.
\begin{theorem}[Zeilberger]
\label{thmASM}
The number of $ n\times n $ alternating sign matrices is
\begin{equation}
\label{eqASM}
    A(n) = \prod_{j=0}^{n-1}\frac{(3j+1)!}{(n+j)!}.
\end{equation}
\end{theorem}
The first proof was given by Zeilberger \cite[p.~5]{Zei96}. However, shortly afterward a simpler proof was found by Kuperberg \cite[Theorem~1]{Kup96}, which exploits a bijection to a ``solvable lattice model'' studied in statistical mechanics.

\begin{figure}[ht]
\centering
\begin{tikzpicture}[scale=0.75]
	\draw[-] (0,3)--(1,3) node[pos=0]{\arrowRight};
	\draw[-] (1,3)--(2,3) node[pos=0.5]{\arrowLeft};
	\draw[-] (2,3)--(3,3) node[pos=0.5]{\arrowLeft};
	\draw[-] (3,3)--(4,3) node[pos=0.5]{\arrowLeft};
	\draw[-] (4,3)--(5,3) node[pos=0.5]{\arrowLeft};
	\draw[-] (5,3)--(6,3) node[pos=1]{\arrowLeft};
	\draw[-] (0,2)--(1,2) node[pos=0]{\arrowRight};
	\draw[-] (1,2)--(2,2) node[pos=0.5]{\arrowRight};
	\draw[-] (2,2)--(3,2) node[pos=0.5]{\arrowLeft};
	\draw[-] (3,2)--(4,2) node[pos=0.5]{\arrowRight};
	\draw[-] (4,2)--(5,2) node[pos=0.5]{\arrowRight};
	\draw[-] (5,2)--(6,2) node[pos=1]{\arrowLeft};
	\draw[-] (0,1)--(1,1) node[pos=0]{\arrowRight};
	\draw[-] (1,1)--(2,1) node[pos=0.5]{\arrowRight};
	\draw[-] (2,1)--(3,1) node[pos=0.5]{\arrowRight};
	\draw[-] (3,1)--(4,1) node[pos=0.5]{\arrowLeft};
	\draw[-] (4,1)--(5,1) node[pos=0.5]{\arrowLeft};
	\draw[-] (5,1)--(6,1) node[pos=1]{\arrowLeft};
	\draw[-] (1,0)--(1,1) node[pos=0]{\arrowDown};
	\draw[-] (1,1)--(1,2) node[pos=0.5]{\arrowDown};
	\draw[-] (1,2)--(1,3) node[pos=0.5]{\arrowDown};
	\draw[-] (1,3)--(1,4) node[pos=1]{\arrowUp};
	\draw[-] (2,0)--(2,1) node[pos=0]{\arrowDown};
	\draw[-] (2,1)--(2,2) node[pos=0.5]{\arrowDown};
	\draw[-] (2,2)--(2,3) node[pos=0.5]{\arrowUp};
	\draw[-] (2,3)--(2,4) node[pos=1]{\arrowUp};
	\draw[-] (3,0)--(3,1) node[pos=0]{\arrowDown};
	\draw[-] (3,1)--(3,2) node[pos=0.5]{\arrowUp};
	\draw[-] (3,2)--(3,3) node[pos=0.5]{\arrowDown};
	\draw[-] (3,3)--(3,4) node[pos=1]{\arrowDownRed};
	\draw[-] (4,0)--(4,1) node[pos=0]{\arrowDown};
	\draw[-] (4,1)--(4,2) node[pos=0.5]{\arrowDown};
	\draw[-] (4,2)--(4,3) node[pos=0.5]{\arrowDown};
	\draw[-] (4,3)--(4,4) node[pos=1]{\arrowDownRed};
	\draw[-] (5,0)--(5,1) node[pos=0]{\arrowDown};
	\draw[-] (5,1)--(5,2) node[pos=0.5]{\arrowDown};
	\draw[-] (5,2)--(5,3) node[pos=0.5]{\arrowUp};
	\draw[-] (5,3)--(5,4) node[pos=1]{\arrowUp};
\end{tikzpicture}
\caption{One possible state of a six-vertex model on the $ 3\times 5 $ grid.}
\label{fig1}
\end{figure}
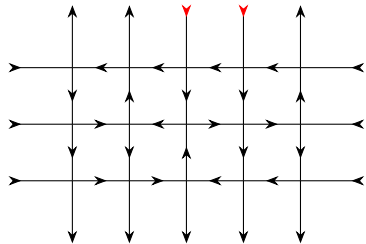

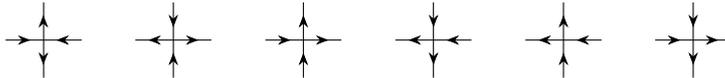
\begin{figure}[ht]
\centering
\begin{subfigure}[b]{0.1\textwidth}
\begin{tikzpicture}[scale=0.5]
	\draw[-] (0,1)--(1,1) node[pos=0.5]{\arrowRight};
	\draw[-] (1,1)--(2,1) node[pos=0.5]{\arrowLeft};
	\draw[-] (1,0)--(1,1) node[pos=0.5]{\arrowDown};
	\draw[-] (1,1)--(1,2) node[pos=0.5]{\arrowUp};
\end{tikzpicture}
\end{subfigure}
\begin{subfigure}[b]{0.1\textwidth}
\begin{tikzpicture}[scale=0.5]
	\draw[-] (0,1)--(1,1) node[pos=0.5]{\arrowLeft};
	\draw[-] (1,1)--(2,1) node[pos=0.5]{\arrowRight};
	\draw[-] (1,0)--(1,1) node[pos=0.5]{\arrowUp};
	\draw[-] (1,1)--(1,2) node[pos=0.5]{\arrowDown};
\end{tikzpicture}
\end{subfigure}
\begin{subfigure}[b]{0.1\textwidth}
\begin{tikzpicture}[scale=0.5]
	\draw[-] (0,1)--(1,1) node[pos=0.5]{\arrowRight};
	\draw[-] (1,1)--(2,1) node[pos=0.5]{\arrowRight};
	\draw[-] (1,0)--(1,1) node[pos=0.5]{\arrowUp};
	\draw[-] (1,1)--(1,2) node[pos=0.5]{\arrowUp};
\end{tikzpicture}
\end{subfigure}
\begin{subfigure}[b]{0.1\textwidth}
\begin{tikzpicture}[scale=0.5]
	\draw[-] (0,1)--(1,1) node[pos=0.5]{\arrowLeft};
	\draw[-] (1,1)--(2,1) node[pos=0.5]{\arrowLeft};
	\draw[-] (1,0)--(1,1) node[pos=0.5]{\arrowDown};
	\draw[-] (1,1)--(1,2) node[pos=0.5]{\arrowDown};
\end{tikzpicture}
\end{subfigure}
\begin{subfigure}[b]{0.1\textwidth}
\begin{tikzpicture}[scale=0.5]
	\draw[-] (0,1)--(1,1) node[pos=0.5]{\arrowLeft};
	\draw[-] (1,1)--(2,1) node[pos=0.5]{\arrowLeft};
	\draw[-] (1,0)--(1,1) node[pos=0.5]{\arrowUp};
	\draw[-] (1,1)--(1,2) node[pos=0.5]{\arrowUp};
\end{tikzpicture}
\end{subfigure}
\begin{subfigure}[b]{0.1\textwidth}
\begin{tikzpicture}[scale=0.5]
	\draw[-] (0,1)--(1,1) node[pos=0.5]{\arrowRight};
	\draw[-] (1,1)--(2,1) node[pos=0.5]{\arrowRight};
	\draw[-] (1,0)--(1,1) node[pos=0.5]{\arrowDown};
	\draw[-] (1,1)--(1,2) node[pos=0.5]{\arrowDown};
\end{tikzpicture}
\end{subfigure}
\caption{The six states of a vertex.}
\label{fig2}
\end{figure}

Roughly speaking, a lattice model is a grid together with all ``allowed'' ways of filling its edges with arrows; particular fillings are called \emph{states of the model} (\Cref{fig1}). In the present paper, we need only consider rectangular grid, together with the rule that every internal vertex where two lines meet must have two arrows entering the vertex, and two arrows leaving. These are usually called \emph{six-vertex models} because there are $\binom{4}{2}=6$ configurations of arrows that satisfy this property. The configuration of arrows adjacent to a vertex is called the \emph{state of the vertex}~(\Cref{fig2}).

The six-vertex model originated in statistical mechanics as a two-dimensional model of ice. One of the simplest cases, known as \emph{square ice}, was solved by Lieb in 1967~\cite{Lie67} by the \emph{Bethe ansatz}, followed by a more general solution by Sutherland~\cite{Sut67}. The energy of the system is captured by the so-called \emph{partition function} (here the word partition is used in a different sense than in the rest of the paper). One cares about the behavior of this function in the \emph{thermodynamic limit}, i.e., as the size of the grid grows to infinity. More information about lattice models from the viewpoint of statistical mechanics can be found in the books~\cite{Bax82, KBI93}.

The number of states of the model depends on its boundary conditions. The six-vertex model considered by Kuperberg~\cite{Kup96}, uses a square lattice and the \emph{domain-wall boundary conditions} (DWBC):
\begin{enumerate}
\item[(a)] The arrows on the left and right boundary point inward.
\item[(b)] The arrows on the bottom boundary point down.
\item[(c)] The arrows on the upper boundary point up.
\end{enumerate}

The states of the six-vertex model with DWBC on the square lattice turn out to be in bijective correspondence with ASMs, and Kuperberg's paper in fact enumerates the former.

Changing the boundary conditions can still yield models of interest. For instance, the \emph{vertically symmetric alternating-sign matrices} (VSASM) of size $(2n+1)\times (2n+1) $ are in bijection with the states of a rectangular $ n \times (2n-1) $ six-vertex model where boundary conditions (a) and (b) remain in place, but (c) has been replaced by:
\begin{enumerate}
    \item[(c')] The arrows on the upper boundary point alternately up and down.
\end{enumerate}

VSASMs are also enumerated by a product formula, conjectured by Mills (see~\cite{Rob91}) and proved by Kuperberg~\cite[Theorem~2]{Kup02}, which appears in~\cite[p.~17]{RS04} in the following form.
\begin{theorem}[Kuperberg]
\label{thmVSASM}
The number of $ (2n+1)\times (2n+1) $ VSASMs is given by
\begin{equation}
\label{eqVSASM}
    A_V(2n+1) = \frac{1}{2^n}\prod_{j=1}^n \frac{(6j-2)!(2j-1)!}{(4j-1)!(4j-2)!}.
\end{equation}
\end{theorem}

In this paper we generalize these results by weakening condition (c) above. The orientation of arrows on the upper boundary will be determined by a partition $ \lambda $. This way of prescribing the boundary condition was used in Brubaker, Bump, and Friedberg~\cite[p.~287]{BBF11} to give a lattice model interpretation for Schur functions.

A \emph{partition} $ \lambda $ of length $ n $ is an $ n $-tuple of integers $ (\lambda_1, \lambda_2, \ldots, \lambda_n) $ such that $ \lambda_1 \geq \lambda_2 \geq \ldots \geq \lambda_n \geq 0 $. Beginning with a fixed $\lambda $ this model uses a rectangular lattice with $n$ rows and $n+\lambda_1$ columns. The left, right, and bottom boundaries are the same as for the DWBC, but the upper depends on $\lambda$ as follows. First, number the columns in increasing order from right to left by the numbers $ 1, 2, \ldots, n+\lambda_1 $. Define $\rho$ to be the partition $(n, n-1, \ldots, 1)$, and $\lambda+\rho$ to be their coordinate-wise sum. Then the arrow on the upper boundary in the $i^\text{th}$ column points up if $i$ is contained in $\lambda+\rho$ and down otherwise.

This has the effect that the arrow in the $ (n+\lambda_1)^{\text{th}} $ (i.e. the leftmost) column points up and the orientation of the remaining arrows on the upper boundary can be prescribed arbitrarily.

\begin{example}
Let $ n = 3 $ and $\lambda = (2, 2, 0) $. Then $ \lambda+\rho = (5, 4, 1) $, yielding the boundary conditions shown in \Cref{fig3}.
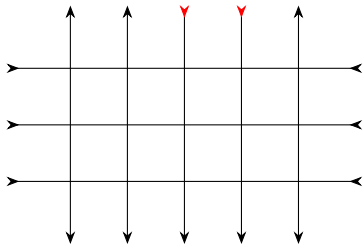
\begin{figure}[ht]
\centering
\begin{tikzpicture}[scale=0.75]
	\draw[-] (0,1)--(6,1) node[pos=0]{\arrowRight} node[pos=1]{\arrowLeft};
	\draw[-] (0,2)--(6,2) node[pos=0]{\arrowRight} node[pos=1]{\arrowLeft};
	\draw[-] (0,3)--(6,3) node[pos=0]{\arrowRight} node[pos=1]{\arrowLeft};
	\draw[-] (1,0)--(1,4) node[pos=0]{\arrowDown} node[pos=1]{\arrowUp};
	\draw[-] (2,0)--(2,4) node[pos=0]{\arrowDown} node[pos=1]{\arrowUp};
	\draw[-] (3,0)--(3,4) node[pos=0]{\arrowDown} node[pos=1]{\arrowDownRed};
	\draw[-] (4,0)--(4,4) node[pos=0]{\arrowDown} node[pos=1]{\arrowDownRed};
	\draw[-] (5,0)--(5,4) node[pos=0]{\arrowDown} node[pos=1]{\arrowUp};
\end{tikzpicture}
\caption{The lattice for $ n = 3 $ and $ \lambda = (2, 2, 0) $.}
\label{fig3}
\end{figure}
\end{example}

\begin{definition} 
Write $A_\lambda(n)$ to denote the number of states of the six-vertex model associated with $\lambda$.
\end{definition}

Our main objective is to prove the following theorem, which gives a broad description of how $A_\lambda$ changes as the first part of $\lambda$ grows:

\begin{theorem}
\label{thmMain}
Let $ \lambda = (\lambda_1, \lambda_2, \ldots, \lambda_n) $ be a partition. Then $ A_\lambda(n) $ is a polynomial in $ \lambda_1 $ of degree $ n-1 $.
\end{theorem}

The proof gives a partial computation for these polynomials, and for two families of $\lambda$ we can complete this computation. We call the partitions in these families \emph{hook shapes} and \emph{staircases} based on the shape of their Ferrers diagrams. For the family of hook shapes, we have the following theorem:

\begin{theorem}
\label{thmHooks}
Let $ m, d \geq 0 $ and let $ \lambda = (m+d, d, \ldots, d) $ be a partition with $n$ parts. Then
\begin{align*}
    A_\lambda(n)& = 
        \frac{A(n)}{\binom{3n-2}{n-1}}\sum_{j=1}^n \binom{m+j-1}{m}\binom{n+j-2}{n-1}\binom{2n-1-j}{n-1}\\
	& = A(n)\sum_{k=0}^{m}\binom{m}{k}\frac{(n+k-1)!(2n-1)!}{k!(n-k-1)!(2n+k-1)!}\\
	& = A(n)\cdot\frac{P_{m}(n)}{Q_{m}(n)}
\end{align*}
where $ P_{m} $ is a polynomial of degree $2m-\lfloor\frac{m+1}{2}\rfloor$, $ Q_{m} $ is a polynomial of degree $m-\lfloor\frac{m+1}{2}\rfloor$ and $ A(n) $ is given by~\Cref{eqASM}.
\end{theorem}

The first sum in this theorem can be used, for a fixed $ n $, to explicitly compute $ A_\lambda(n) $ as a polynomial in $ m $. Thus it describes what happens when we increase the arm length of the hook. The second sum can be used, for a fixed $ m $, to compute $ A_\lambda(n) $ as a rational function in $ n $, and thus it describes what happens as the leg length increases.

In addition to hook shapes, we find also a similar formula for staircases:

\begin{theorem}
\label{thmStairs}
Let $ \lambda = (\lambda_1+d, n-2+d, n-3+d, \ldots, 1+d, d) $ where $ \lambda_1 \geq n-1 $ and $ d \geq 0 $. Then
\[
    A_\lambda(n) = \frac{A_V(2n-1)}{\binom{4n-2}{2n-1}}\sum_{j = 1}^n \binom{\lambda_1+1-j}{\lambda_1+1-n}\binom{2n+j-2}{2n-1}\binom{4n-j-1}{2n-1}
\]
where $ A_V(2n-1) $ is determined by~\Cref{eqVSASM}.
\end{theorem}

We remark that in the preceding theorems, $ A_\lambda(n) $ does not depend on $ d $. In general, as we will observe later, adding an integer $ d \geq 0 $ to all parts of the partition $ \lambda $ does not change $ A_\lambda(n) $.

\

The rest of the paper is organized as follows. In~\Cref{sectionMain} we prove~\Cref{thmMain}. In the course of the proof, we find a general (but not completely explicit) formula for $ A_\lambda(n) $. This formula is then used to prove~\Cref{thmHooks} in~\Cref{sectionHooks} and~\Cref{thmStairs} in~\Cref{sectionStairs}.

\section*{Acknowledgements}
This work originated as a project at the Student Number Theory Seminar at the Charles University. We wish to thank V\'{i}t\v{e}zslav Kala for his guidance and mentorship.

\section{Proof of the Main Theorem}
\label{sectionMain}

We consider first the six-vertex model on an $r\times c$ lattice where $ r, c \geq 1 $, satisfying one of two types of boundary conditions (\Cref{fig4}). Let $ S(r, c) $ be the number of states of a model with the following boundary conditions:
\begin{itemize}
\item The arrows on the left boundary point right.
\item The last arrow on the right boundary points left and all the other ones point right.
\item The arrows on the bottom boundary point down.
\item The first arrow on the upper boundary points up and all the other ones point down.
\end{itemize}
Additionally, let $ T(r, c) $ be the number of states of a model with the following boundary conditions:
\begin{itemize}
\item The first arrow on the left boundary points left and all the other ones point right.
\item The last arrow on the right boundary points left and all the other ones point right.
\item The arrows on the upper and bottom boundaries point down.
\end{itemize}
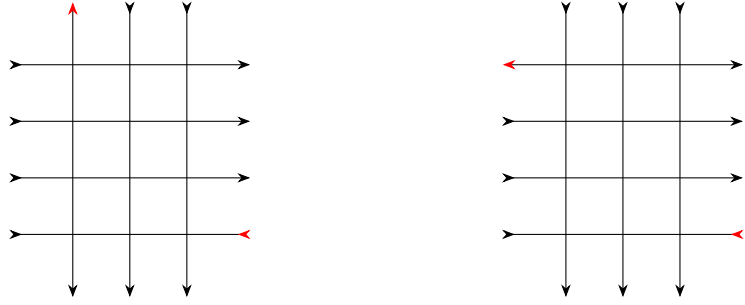
\begin{figure}[ht]
\centering
\begin{subfigure}{0.4\textwidth}
\centering
\begin{tikzpicture}[scale=0.75]
	\draw[-] (0,1)--(4,1) node[pos=0]{\arrowRight} node[pos=1]{\arrowLeftRed};
	\draw[-] (0,2)--(4,2) node[pos=0]{\arrowRight} node[pos=1]{\arrowRight};
	\draw[-] (0,3)--(4,3) node[pos=0]{\arrowRight} node[pos=1]{\arrowRight};
	\draw[-] (0,4)--(4,4) node[pos=0]{\arrowRight} node[pos=1]{\arrowRight};
	\draw[-] (1,0)--(1,5) node[pos=0]{\arrowDown} node[pos=1]{\arrowUpRed};
	\draw[-] (2,0)--(2,5) node[pos=0]{\arrowDown} node[pos=1]{\arrowDown};
	\draw[-] (3,0)--(3,5) node[pos=0]{\arrowDown} node[pos=1]{\arrowDown};
\end{tikzpicture}
\end{subfigure}
\begin{subfigure}{0.4\textwidth}
\centering
\begin{tikzpicture}[scale=0.75]
	\draw[-] (0,1)--(4,1) node[pos=0]{\arrowRight} node[pos=1]{\arrowLeftRed};
	\draw[-] (0,2)--(4,2) node[pos=0]{\arrowRight} node[pos=1]{\arrowRight};
	\draw[-] (0,3)--(4,3) node[pos=0]{\arrowRight} node[pos=1]{\arrowRight};
	\draw[-] (0,4)--(4,4) node[pos=0]{\arrowLeftRed} node[pos=1]{\arrowRight};
	\draw[-] (1,0)--(1,5) node[pos=0]{\arrowDown} node[pos=1]{\arrowDown};
	\draw[-] (2,0)--(2,5) node[pos=0]{\arrowDown} node[pos=1]{\arrowDown};
	\draw[-] (3,0)--(3,5) node[pos=0]{\arrowDown} node[pos=1]{\arrowDown};
\end{tikzpicture}
\end{subfigure}
\caption{The two types of boundary conditions for $ r = 4 $ and $ c = 3 $.}
\label{fig4}
\end{figure}

\begin{lemma}
\label{lemmaS}
If $ r, c \geq 1 $ are integers, then
\[
	S(r, c) = T(r, c) = \binom{r+c-2}{c-1}.
\]
\end{lemma}
\begin{proof}
If $ r = 1 $ or $ c = 1 $, one easily observes that there is a unique state:
\[
	S(1, c) = S(r, 1) = T(1, c) = T(r, 1) = 1.
\]
From now on we assume $r, c>1$ and we find a recurrence for $ S(r, c) $ and $ T(r, c) $.

Take the model with the first type of boundary conditions. In the first column from the left, we can choose from two possible states for the first vertex. In the first case, the arrows point up and right. This forces the orientation in the entire first row and the number of states we get this way is $ S(r-1, c) $. In the second case, the arrows point down and left. This forces the orientation in the first column and the number of states obtained this way is $ T(r, c-1) $. Thus
\[
	S(r, c) = S(r-1, c)+T(r, c-1).
\]

Now we apply similar reasoning to the model with the second type of boundary conditions. Again, we can choose from two possible states for the first vertex in the first column on the left. Either the arrows point up and right, in which case we get $ S(r-1, c) $ states. Or the arrows point down and left, and we get $ T(r, c-1) $ states. Thus
\[
	T(r, c) = S(r-1, c)+T(r, c-1).
\]

From the two recurrences together with the equalities $ S(1, c) = T(1, c) $ and $ S(r, 1) = T(r, 1) $ above, it follows by induction that $ S(r, c) = T(r, c) $. To summarize, we proved
\begin{align*}
	S(r, c)& = S(r-1, c)+S(r, c-1),\qquad r, c \geq 2\\
	S(r, 1)& = S(1, c) = 1.
\end{align*}

The recurrence is Pascal's rule and the initial conditions agree with $ \binom{r+c-2}{c-1} $.
\end{proof}

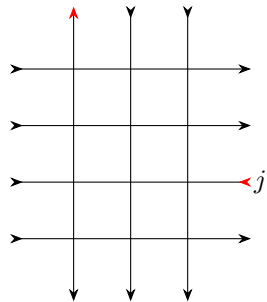
\begin{figure}[hb]
\centering
\begin{tikzpicture}[scale=0.75]
	\draw[-] (0,4)--(4,4) node[pos=0]{\arrowRight} node[pos=1]{\arrowRight};
	\draw[-] (0,3)--(4,3) node[pos=0]{\arrowRight} node[pos=1]{\arrowRight};
	\draw[-] (0,2)--(4,2) node[pos=0]{\arrowRight} node[pos=1]{\arrowLeftRed};
	\draw[-] (0,1)--(4,1) node[pos=0]{\arrowRight} node[pos=1]{\arrowRight};
	\draw[-] (1,0)--(1,5) node[pos=0]{\arrowDown} node[pos=1]{\arrowUpRed};
	\draw[-] (2,0)--(2,5) node[pos=0]{\arrowDown} node[pos=1]{\arrowDown};
	\draw[-] (3,0)--(3,5) node[pos=0]{\arrowDown} node[pos=1]{\arrowDown};
	\draw[-] node[right] at (4,2) {$j$};
\end{tikzpicture}
\caption{The $ n \times (m+1) $ lattice with a unique arrow pointing left on the right boundary in the $ j^\text{th} $ row.}
\label{fig5}
\end{figure}

Next, let $ n \geq 1 $, $ m \geq 0 $, and $ 1 \leq j \leq n $. Consider a lattice with $ n $ rows and $ m+1 $ columns with the following boundary conditions (\Cref{fig5}):
\begin{itemize}
\item The first arrow on the upper boundary points up and all the other ones point down.
\item The arrows on the bottom boundary point down.
\item The arrows on the left boundary point right.
\item There is a unique arrow pointing left on the right boundary in the $ j^\text{th} $ row; all the other arrows on the right boundary point right.
\end{itemize}

\begin{definition}
Write $ L(m, j) $ for the number of states of this six-vertex model.
\end{definition}

\begin{lemma}
\label{lemmaS1}
If $ n \geq 1 $, $ m \geq 0 $, and $ 1 \leq j \leq n $, then
\[
	L(m, j) = \binom{m+j-1}{m}.
\]
\end{lemma}
\begin{proof}
The states of vertices below the $ j^\text{th} $ row are uniquely determined, and all the vertical arrows adjacent to this row point down. The remaining portion is the $r \times c$ lattice considered above, for $r=j$ and $c=m+1$, which yields the desired formula
\[
	L(m, j) = S(j, m+1) = \binom{m+j-1}{m}.\qedhere
\]
\end{proof}

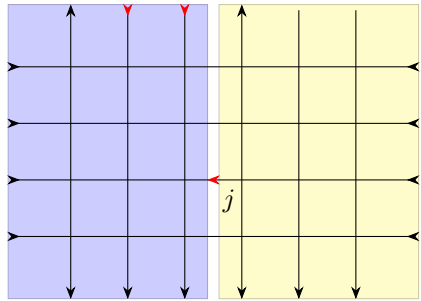
\begin{figure}[ht]
\centering
\begin{tikzpicture}[scale=0.75]
	\draw[fill=blue, opacity=0.2] (-0.1,-0.1) rectangle (3.4,5.1);
	\draw[fill=yellow, opacity=0.2] (3.6, -0.1) rectangle (7.1, 5.1);
	\draw[-] (0,1)--(7,1) node[pos=0]{\arrowRight} node[pos=1]{\arrowLeft};
	\draw[-] (0,2)--(3,2) node[pos=0]{\arrowRight};
	\draw[-] (3,2)--(4,2) node[pos=0.5]{\arrowLeftRed};
	\draw[-] (4,2)--(7,2) node[pos=1]{\arrowLeft};
	\draw[-] (0,3)--(7,3) node[pos=0]{\arrowRight} node[pos=1]{\arrowLeft};
	\draw[-] (0,4)--(7,4) node[pos=0]{\arrowRight} node[pos=1]{\arrowLeft};
	\node [below right] at (3.5, 2) {$j$};
	\draw[-] (1,0)--(1,5) node[pos=0]{\arrowDown} node[pos=1]{\arrowUp};
	\draw[-] (2,0)--(2,5) node[pos=0]{\arrowDown} node[pos=1]{\arrowDownRed};
	\draw[-] (3,0)--(3,5) node[pos=0]{\arrowDown} node[pos=1]{\arrowDownRed};
	\draw[-] (4,0)--(4,5) node[pos=0]{\arrowDown} node[pos=1]{\arrowUp};
	\draw[-] (5,0)--(5,5) node[pos=0]{\arrowDown};
	\draw[-] (6,0)--(6,5) node[pos=0]{\arrowDown};
\end{tikzpicture}
\caption{The division of the lattice into two parts.}
\label{fig6}
\end{figure}

Let $\lambda = (\lambda_1, \lambda_2, \ldots, \lambda_n) $ be a partition of length $n$. We divide the lattice with $n$ rows and $n+\lambda_1$ columns into two parts (see~\Cref{fig6}):
\begin{enumerate}
    \item[(L)] The left part consists of $ 1+\lambda_1-\lambda_2$ columns.
    \item[(R)] The right part contains the remaining $ n+\lambda_2-1 $ columns.
\end{enumerate}

On the upper boundary of the left part, the arrow in the $(n+\lambda_1)^\text{th}$ column points up and the arrows in the remaining $ \lambda_1-\lambda_2 $ columns point down. In each state, the left part has exactly one arrow pointing left on its right boundary. As before, let $ 1 \leq j \leq n $ and assume that the unique arrow pointing left on the right boundary of the left part is in the $ j^\text{th} $ row.

The number of states of the left part is given by $ L(\lambda_1-\lambda_2, j) $. Additionally, let us define $ R(\lambda, j) $ to be the number of states of the right part. Note that while $R$ depends upon $\lambda$, it is independent of $\lambda_1$.

\begin{lemma}
\label{lemmaMain}
Let $ n \geq 2 $ and let $ \lambda = (\lambda_1, \lambda_2, \ldots, \lambda_n) $ be a partition with $ n $ parts. Then
\[
	A_\lambda(n) = \sum_{j=1}^n\binom{\lambda_1-\lambda_2+j-1}{\lambda_1-\lambda_2}R(\lambda, j).
\]
\end{lemma}
\begin{proof}
Let $ m = \lambda_1-\lambda_2 $. By running through all possible choices for $ j $, we get
\[
	A_\lambda(n) = \sum_{j=1}^n L(m, j)R(\lambda, j).
\]
But we know from Lemma~\ref{lemmaS1} that 
\[
	L(m, j) = \binom{m+j-1}{m}.
\]
Thus
\[
	A_\lambda(n) = \sum_{j=1}^n\binom{m+j-1}{m}R(\lambda, j).\qedhere
\]
\end{proof}

This lemma completes the proof of~\Cref{thmMain}, since $ \binom{\lambda_1-\lambda_2+j-1}{\lambda_1-\lambda_2} $ is a polynomial in $ \lambda_1 $ of degree $ j-1 $. Thus the sum has degree $ n-1 $ as desired.

\section{Hook shapes}
\label{sectionHooks}

\begin{figure}[ht]
\centering
\begin{tikzpicture}[scale=0.75]
	\draw[-] (0,3)--(5,3) node[pos=0]{\arrowRight};
	\draw[-] (5,3)--(6,3) node[pos=0.5]{\arrowLeft};
	\draw[-] (6,3)--(7,3) node[pos=1]{\arrowLeft};
	\draw[-] (0,2)--(5,2) node[pos=0]{\arrowRight};
	\draw[-] (5,2)--(6,2) node[pos=0.5]{\arrowLeft};
	\draw[-] (6,2)--(7,2) node[pos=1]{\arrowLeft};
	\draw[-] (0,1)--(5,1) node[pos=0]{\arrowRight};
	\draw[-] (5,1)--(6,1) node[pos=0.5]{\arrowLeft};
	\draw[-] (6,1)--(7,1) node[pos=1]{\arrowLeft};
	\draw[-] (1,4)--(1,0) node[pos=0]{\arrowUp} node[pos=1]{\arrowDown};
	\draw[-] (2,4)--(2,0) node[pos=0]{\arrowDown} node[pos=1]{\arrowDown};
	\draw[-] (3,4)--(3,0) node[pos=0]{\arrowDown} node[pos=1]{\arrowDown};
	\draw[-] (4,4)--(4,0) node[pos=0]{\arrowUp} node[pos=1]{\arrowDown};
	\draw[-] (5,4)--(5,0) node[pos=0]{\arrowUp} node[pos=1]{\arrowDown};
	\draw[-] (6,4)--(6,3) node[pos=0]{\arrowDown};
	\draw[-] (6,3)--(6,2) node[pos=0.5]{\arrowDown};
	\draw[-] (6,2)--(6,1) node[pos=0.5]{\arrowDown};
	\draw[-] (6,1)--(6,0) node[pos=1]{\arrowDown};
\end{tikzpicture}
\caption{A $3\times 6 $ lattice corresponding to the hook shape $ (3, 1, 1) $. Its $ 3\times 5 $ sublattice corresponds to the partition $ (2, 0, 0) $.}
\label{fig7}
\end{figure}
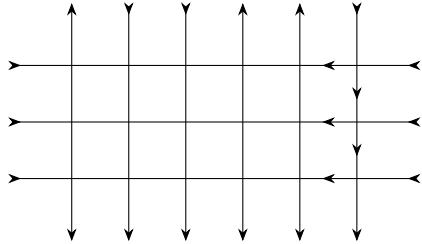

Let us first make an easy observation: If $ \lambda = (\lambda_1, \lambda_2, \ldots, \lambda_n) $ and $ \vec{1} = (1, 1, \ldots, 1) $ are partitions with $ n $ parts, then
\[
	A_\lambda(n) = A_{\lambda+\vec{1}}(n).
\]
Indeed, adding 1 to all the entries simply appends a column on the right with arrows pointing down, and together with the right boundary condition this implies that that column is forced~(\Cref{fig7}).

We turn to the proof of~\Cref{thmHooks}. By the preceding observation, we can assume $ d = 0 $ and work with the partition $ \lambda = (m, 0, \ldots, 0) $. 

Let $ n \geq 2 $ and $ 1 \leq j \leq n $. Consider an $ (n-1)\times n $ lattice satisfying the following boundary conditions:
\begin{itemize}
\item The arrows on the left and right boundary point inward.
\item The arrows on the bottom boundary point down.
\item The arrow in the $ j^\text{th} $ column on the upper boundary points down; the remaining arrows on the upper boundary point up.
\end{itemize}

The states of this six-vertex model are in bijection with $ n \times n $ ASMs where the unique 1 in the first row is in the $ j^\text{th} $ column. Fortunately, the number of these matrices has already been computed. This is precisely the \emph{Refined ASM Theorem} proved by Zeilberger~\cite[p.~60]{Zei96a}.
\begin{theorem}[Zeilberger]
\label{thmRefinedASM}
Let $ 1 \leq j \leq n $. There are
\[
A(n, j) = A(n)\frac{\binom{n+j-2}{n-1}\binom{2n-1-j}{n-1}}{\binom{3n-2}{n-1}}.
\]
$ n\times n $ alternating sign matrices such that the unique 1 in the first row is in the $ j^\text{th} $ column.
\end{theorem}

The theorem is used in the proof of the next lemma.

\begin{lemma}
\label{lemmaHooks1}
If $ \lambda = (m, 0, \ldots, 0) $ is a partition with $ n $ parts, then
\[
	A_\lambda(n) = \frac{A(n)}{\binom{3n-2}{n-1}}\sum_{j=1}^n \binom{m+j-1}{m}\binom{n+j-2}{n-1}\binom{2n-1-j}{n-1}.
\]
\end{lemma}
\begin{proof}
We can assume $ n \geq 2 $ because if $ n = 1 $, then $ A_\lambda(1) = 1 $ and the right side of the formula also equals $1$.

The lattice has dimensions $ n \times (n+m) $. Setting $ \rho = (n, n-1, \ldots, 1) $, we have $ \lambda+\rho = (n+m, n-1, \ldots, 1) $. The orientation of arrows on the upper boundary (counted from the right) is as follows: the arrows in columns $1, 2, \ldots, n-1$ point up, the arrows in the next $ m $ columns point down and the arrow in column $ (n+m) $ points up. As before, we divide the lattice into two parts where the left part consists of $ m+1 $ columns and the right part of $ n-1 $ columns.

Now we consider the right part. We know that there is a unique arrow pointing left on its left boundary. Let this arrow be placed in the $ j^\text{th} $ row where $ 1 \leq j \leq n $. Rotate clockwise by 90 degrees and reverse the directions of all the arrows; the resulting rectangular lattice of size $(n-1)\times n$ has domain wall boundary conditions, except the arrow in the $j^\text{th}$ column (counted from the right) on the upper boundary has been flipped. Thus $ R(\lambda, j) = A(n, j) $.

By~\Cref{lemmaMain} and~\Cref{thmRefinedASM}, the total number of states is
\[
	A_\lambda(n) = \sum_{j=1}^n \binom{m+j-1}{m}A(n,j) = \frac{A(n)}{\binom{3n-2}{n-1}}\sum_{j=1}^n \binom{m+j-1}{m}\binom{n+j-2}{n-1}\binom{2n-1-j}{n-1}.\qedhere
\]
\end{proof}

In the following, we adopt the usual conventions regarding binomial coefficients. For integers $ n $ and $ k $, we have $ \binom{n}{k} = 1 $ if $ k = 0 $ and $ \binom{n}{k} = 0 $ if $ k < 0 $. We also make a frequent use of the generalized binomial coefficient defined as
\[
	\binom{\alpha}{k} = \frac{\alpha(\alpha-1)\ldots (\alpha-k+1)}{k!}
\]
for any complex number $ \alpha $ and integer $ k \geq 0 $.

Let $ \lambda = (m, 0, \ldots, 0) $ be a partition with $ n $ parts. Shifting the index of summation, we write the formula in~\Cref{lemmaHooks1} as
\[
	A_\lambda(n) = \frac{A(n)}{\binom{3n-2}{n-1}}\sum_{k=0}^{n-1} \binom{m+k}{m}\binom{n+k-1}{n-1}\binom{2n-k-2}{n-1}.
\]

The last two binomial coefficients in the formula for $ A_\lambda(n) $ can be expressed as
\begin{align*}
	\binom{n+k-1}{n-1}& = \frac{n(n+1)\ldots(n+k-1)}{k!} = (-1)^k\frac{(-n)(-n-1)\ldots(-n-k+1)}{k!} = (-1)^k\binom{-n}{k},\\
	\binom{2n-k-2}{n-1}& = \frac{n(n+1)\ldots(2n-k-2)}{(n-k-1)!} = (-1)^{n-k-1}\frac{(-n)(-n-1)\ldots(-n-(n-k-2))}{(n-k-1)!}\\
	& = (-1)^{n-1-k}\binom{-n}{n-1-k},
\end{align*}
which allows us to write
\begin{equation}
\label{eqAlambdaInter}
	A_\lambda(n) = \frac{(-1)^{n-1}A(n)}{\binom{3n-2}{n-1}}\sum_{k=0}^{n-1}\binom{m+k}{m}\binom{-n}{k}\binom{-n}{n-1-k}.
\end{equation}

The last sum satisfies the following combinatorial identity. The proof is a standard application of generating functions.
\begin{lemma}
\label{lemmaIdentity}
If $ m \geq 0 $ and $ n \geq 1 $ are integers, then
\[
	\sum_{k=0}^{n-1} \binom{m+k}{m}\binom{-n}{k}\binom{-n}{n-1-k} = \sum_{k=0}^{m}\binom{m}{k}\binom{-n}{k}\binom{-2n-k}{n-1-k}.
\]
\end{lemma}
\begin{proof}
Consider the function
\[
	\dv[m]{x}(\frac{x^m}{m!}(1+x)^{-n}) x^{1-n}(1+x)^{-n}.
\]
We write it in two ways. First,
\begin{align*}
	\dv[m]{x}(\frac{x^m}{m!}(1+x)^{-n})& = \dv[m]{x}(\sum_{k=0}^\infty \frac{1}{m!}\binom{-n}{k} x^{m+k})\\
	& = \sum_{k=0}^\infty \frac{(m+k)(m+k-1)\ldots (k+1)}{m!}\binom{-n}{k}x^k\\
	& = \sum_{k=0}^\infty \binom{m+k}{m}\binom{-n}{k}x^k
\end{align*}
and
\[
	x^{1-n}(1+x)^{-n} = \sum_{j=0}^\infty \binom{-n}{j}x^{1-n+j}.
\]
The product of the two series equals
\[
	\sum_{k=0}^\infty \binom{m+k}{m}\binom{-n}{k}x^k\cdot \sum_{j=0}^\infty \binom{-n}{j}x^{1-n+j} = \sum_{i = 0}^\infty \sum_{k = 0}^i\binom{m+k}{m}\binom{-n}{k}\binom{-n}{i-k} x^{1-n+i}.
\]
In particular, the coefficient of the constant term (for $ i = n-1 $) is
\begin{equation}
\label{eqExp1}
	\sum_{k=0}^{n-1} \binom{m+k}{m}\binom{-n}{k}\binom{-n}{n-1-k}.
\end{equation}

Secondly, we take the $ m^\text{th} $ derivative using the Leibniz rule:
\begin{align*}
	\dv[m]{x}(\frac{x^m}{m!}(1+x)^{-n})& = \sum_{k = 0}^m \frac{1}{m!}\binom{m}{k}\dv[m-k]{x^m}{x}\dv[k]{(1+x)^{-n}}{x}\\
	& = \sum_{k = 0}^m \frac{1}{m!}\binom{m}{k}\frac{m!}{k!}x^k(-n)(-n-1)(-n-(k-1))\cdot(1+x)^{-n-k}\\
	& = \sum_{k=0}^m \binom{m}{k}\binom{-n}{k}x^k (1+x)^{-n-k}.
\end{align*}
Hence
\begin{align*}
	\dv[m]{x}(\frac{x^m}{m!}(1+x)^{-n}) x^{1-n}(1+x)^{-n}& = \sum_{k=0}^m \binom{m}{k}\binom{-n}{k}x^{1-n+k} (1+x)^{-2n-k}\\
	& = \sum_{k=0}^m \binom{m}{k}\binom{-n}{k}x^{1-n+k} \sum_{j=0}^\infty \binom{-2n-k}{j} x^j\\
	& = \sum_{i=0}^\infty \sum_{k=0}^m \binom{m}{k}\binom{-n}{k}\binom{-2n-k}{i-k}x^{1-n+i}
\end{align*}
and the constant coefficient (for $ i = n-1 $) is
\begin{equation}
\label{eqExp2}
	\sum_{k=0}^m \binom{m}{k}\binom{-n}{k}\binom{-2n-k}{n-1-k}.
\end{equation}
Comparing~\Cref{eqExp1,eqExp2}, we get the formula.
\end{proof}

\begin{lemma}
\label{lemmaHooks2}
If $ \lambda = (m, 0, \ldots, 0) $ is a partition with $ n $ parts, then
\[
	A_\lambda(n) = A(n)\sum_{k=0}^{m}\binom{m}{k}\frac{1}{k!}\frac{(n-k)(n-k+1)\ldots(n+k-1)}{(2n)(2n+1)\ldots(2n+k-1)}.
\]
\end{lemma}
\begin{proof}
From~\Cref{eqAlambdaInter} and~\Cref{lemmaIdentity}
\[
	A_\lambda(n) = \frac{(-1)^{n-1} A(n)}{\binom{3n-2}{n-1}}\sum_{k=0}^m\binom{m}{k}\binom{-n}{k}\binom{-2n-k}{n-1-k}.
\]
The last step is to show that the right side equals $ A(n) $ times a rational function in $ n $. We have
\begin{align*}
	\frac{(-1)^{n-1}\binom{-n}{k}\binom{-2n-k}{n-1-k}}{\binom{3n-2}{n-1}}& = (-1)^{n-1}\frac{(-n)(-n-1)\ldots(-n-k+1)}{k!}\\
	&\cdot \frac{(-2n-k)(-2n-k-1)\ldots(-3n+2)}{(n-1-k)!}\cdot \frac{(n-1)!}{(2n)(2n+1)\ldots(3n-2)}\\
	& = \frac{1}{k!}\frac{(n-k)(n-k+1)\ldots(n+k-1)}{(2n)(2n+1)\ldots(2n+k-1)}.
\end{align*}
The formula follows.
\end{proof}

We let $ R_m(n) $ denote the rational function factor, so that $ A_\lambda(n) = A(n)R_m(n) $. These factors for $ m \leq 5 $ are given in~\Cref{table1}.

\begin{proof}[Proof of~\Cref{thmHooks}]
By the observation at the beginning of this section, we can assume $ d = 0 $. The first and second equality in the theorem follow from~\Cref{lemmaHooks1} and~\Cref{lemmaHooks2}, respectively.

After factoring $ A(n) $ from the right side, the remaining factor is
\[
	R_m(n) = \sum_{k=0}^m\binom{m}{k}\frac{1}{k!}\frac{(n-k)(n-k+1)\dots(n+k-1)}{(2n)(2n+1)\dots (2n+k-1)}.
\]
Rewriting these rational functions to have a common denominator, we get
\[
	R_m(n) = \frac{\sum_{k=0}^m\binom{m}{k}\frac{1}{k!}(n-k)(n-k+1)\dots(n+k-1)\cdot(2n+k)\dots (2n+m-1)}{(2n)(2n+1)\dots (2n+m-1)}.
\]
Consider the even factors $ (2n+2i) $ where $ 0 \leq i \leq \lfloor\frac{m-1}{2}\rfloor $ in the denominator. For $ 0 \leq i \leq k-1 $, this factor cancels with $ (n+i) $ and for $ k \leq i \leq \lfloor\frac{m-1}{2}\rfloor $, it cancels with $ (2n+2i) $. Thus
\[
	R_m(n) = \frac{P_m(n)}{Q_m(n)}
\]
where $ P_m(n) $ is a polynomial of degree $ 2m-\lfloor\frac{m+1}{2}\rfloor $ and $ Q_m(n) $ is a polynomial of degree $ m-\lfloor\frac{m+1}{2}\rfloor $.
\end{proof}

\begin{table}[ht]
\centering
\begin{TAB}(r,1cm,1.2cm)[5pt]{|c|l|}{|c|c|c|c|c|c|c|}
	$m$&$R_m(n)$\\
	0&1\\
	1&$\frac{n+1}{2}$\\
	2&$\frac{n^3+6n^2+3n+2}{4(2n+1)}$\\
	3&$\frac{n^4+14n^3+35n^2+10n+12}{24(2n+1)}$\\
	4&$\frac{n^6+27n^5+199n^4+456n^3+448n^2+156n+144}{96(4n^2+8n+3)}$\\
	5&$\frac{n^7+42n^6+542n^5+2540n^4+4569n^3+4138n^2+1128n+1440}{960(4n^2+8n+3)}$
\end{TAB}
\caption{The factors $ R_m(n) $ for $ m \leq 5 $.}
\label{table1}
\end{table}

\section{Staircases}
\label{sectionStairs}

In the hook shapes case, we saw that $R(\lambda, j)$ was given by the refined ASM enumeration. As it happens, there is a similar result in the staircase setting.

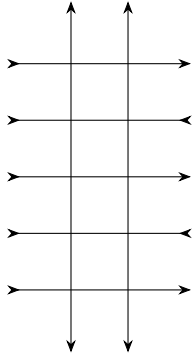
\begin{figure}[ht]
\centering
\begin{tikzpicture}[scale=0.75]
	\draw[-] (0,1)--(3,1) node[pos=0]{\arrowRight} node[pos=1]{\arrowRight};
	\draw[-] (0,2)--(3,2) node[pos=0]{\arrowRight} node[pos=1]{\arrowLeft};
	\draw[-] (0,3)--(3,3) node[pos=0]{\arrowRight} node[pos=1]{\arrowRight};
	\draw[-] (0,4)--(3,4) node[pos=0]{\arrowRight} node[pos=1]{\arrowLeft};
	\draw[-] (0,5)--(3,5) node[pos=0]{\arrowRight} node[pos=1]{\arrowRight};
	\draw[-] (1,0)--(1,6) node[pos=0]{\arrowDown} node[pos=1]{\arrowUp};
	\draw[-] (2,0)--(2,6) node[pos=0]{\arrowDown} node[pos=1]{\arrowUp};
\end{tikzpicture}
\caption{An $ (2n+1)\times n $ lattice whose states are in bijection with $ (2n+1)\times (2n+1) $ VSASMs.} 
\label{fig8}
\end{figure}

This result concerns the \emph{vertically symmetric} alternating-sign matrices (VSASMs). As their name suggests, these are ASMs which are preserved when reflecting entries across the middle column, such as  
\[
	\begin{pmatrix}
		0&0&1&0&0\\
		1&0&-1&0&1\\
		0&0&1&0&0\\
		0&1&-1&1&0\\
		0&0&1&0&0
	\end{pmatrix}.
\]
It is easy to check that the middle column always contains an alternating sequence of $ 1 $'s and $ -1 $'s, and so it follows that $ (2n+1)\times (2n+1) $ VSASMs are in bijection with $ (2n+1) \times n $ lattices of the form indicated in~\Cref{fig8} representing the left half of the matrix.

Observe that the unique $ 1 $ in the first and last row is in the middle, which corresponds to the fact that the orientation of arrows in the first and last row is fixed. Deleting these rows, we get a $ (2n-1)\times n $ lattice with arrows on the right boundary pointing alternately left and right.

As recently as 2021, Fischer and Saikia~\cite[Theorem~3.2]{FS21} published a formula for the refined enumeration of VSASMs conjectured by Fischer~\cite[p.~538]{Fis09}.
\begin{theorem}[Fischer and Saikia]
\label{thmRefinedVSASM}
Let $ 1 \leq i \leq n $. The number of $ (2n+1)\times (2n+1) $ VSASMs with the first 1 in the second row in the $ i^\text{th} $ column is given by
\[
	A_V(2n+1, i) = \frac{\binom{2n+i-2}{2n-1}\binom{4n-i-1}{2n-1}}{\binom{4n-2}{2n-1}}A_V(2n-1).
\]
\end{theorem}

These matrices are in bijection with the states of a six-vertex model on a $ (2n-2)\times n $ lattice satisfying the following boundary conditions:
\begin{itemize}
\item The arrows on the left boundary point right.
\item The arrows on the right boundary point alternately left and right.
\item The arrows on the bottom boundary point down.
\item The arrow in the $ i^\text{th} $ column on the upper boundary points down; the remaining arrows on the upper boundary point up.
\end{itemize}

\begin{lemma}
\label{lemmaStairs}
If $ \lambda = (n-1+m, n-2, \ldots, 1, 0) $ where $ m \geq 0 $, then
\[
	A_\lambda(n) = \sum_{i=1}^n\binom{m+i-1}{m}A_V(2n+1, i).
\]
\end{lemma}
\begin{proof}
We can assume $ n \geq 2 $ because if $ n = 1 $, then $ A_\lambda(1) = 1 $ and $ A_V(3, 1) = 1 $.

Setting $ \rho = (n, n-1, \ldots, 1) $, we have $ \lambda+\rho = (2n-1+m, 2n-3, \ldots, 3, 1) $. The arrows on the upper boundary (counted from the right) are oriented as follows: the arrows in the first $ 2n-3 $ columns point alternately up and down, the arrows in the next $ m+1 $ columns point down, and the arrow in the leftmost column points up.

We again split the lattice into two parts where the left part consists of $ m+1 $ columns and the right part of the remaining $ 2n-2 $ columns. The right part has a unique arrow pointing left on its left boundary (\Cref{fig9}).

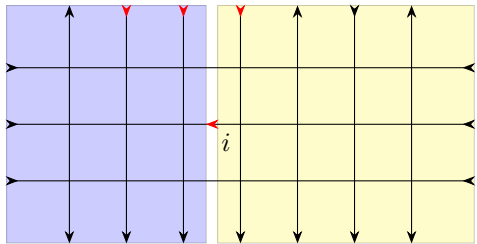
\begin{figure}[ht]
\centering
\begin{tikzpicture}[scale=0.75]
	\draw[fill=blue, opacity=0.2] (-0.1,-0.1) rectangle (3.4,4.1);
	\draw[fill=yellow, opacity=0.2] (3.6, -0.1) rectangle (8.1, 4.1);
	\draw[-] (0,1)--(8,1) node[pos=0]{\arrowRight} node[pos=1]{\arrowLeft};
	\draw[-] (0,2)--(3,2) node[pos=0]{\arrowRight};
	\draw[-] (3,2)--(4,2) node[pos=0.5]{\arrowLeftRed};
	\draw[-] (4,2)--(8,2) node[pos=1]{\arrowLeft};
	\draw[-] (0,3)--(8,3) node[pos=0]{\arrowRight} node[pos=1]{\arrowLeft};
	\node [below right] at (3.5, 2) {$i$};
	\draw[-] (1,0)--(1,4) node[pos=0]{\arrowDown} node[pos=1]{\arrowUp};
	\draw[-] (2,0)--(2,4) node[pos=0]{\arrowDown} node[pos=1]{\arrowDownRed};
	\draw[-] (3,0)--(3,4) node[pos=0]{\arrowDown} node[pos=1]{\arrowDownRed};
	\draw[-] (4,0)--(4,4) node[pos=0]{\arrowDown} node[pos=1]{\arrowDownRed};
	\draw[-] (5,0)--(5,4) node[pos=0]{\arrowDown} node[pos=1]{\arrowUp};
	\draw[-] (6,0)--(6,4) node[pos=0]{\arrowDown} node[pos=1]{\arrowDown};
	\draw[-] (7,0)--(7,4) node[pos=0]{\arrowDown} node[pos=1]{\arrowUp};
\end{tikzpicture}
\caption{The division of the lattice for $ n = 3 $ and $ m = 2 $.}
\label{fig9}
\end{figure}

Let $ 1 \leq i \leq n $ be fixed and suppose that the unique arrow pointing left on the boundary of the right part is in the $ i^\text{th} $ row. The number of states of the left part is given by $ L(m, i) $. Let $ R(\lambda, i) $ be the number of states of the right part. Rotating the right part clockwise by 90 degrees and reversing the orientation of all arrows, we get a $ (2n-2)\times n $ lattice where the unique arrow pointing down on the upper boundary is in the $ i^\text{th} $ column (counted from the right). We know from the preceding discussion that the number of states of this lattice is $ R(\lambda, i) = A_V(2n+1, i) $.

The number of states is given by
\[
	A_\lambda(n) = \sum_{i=1}^n L(m, i)R(\lambda, i) = \sum_{i=1}^n \binom{m+i-1}{m} A_V(2n+1, i).\qedhere
\]
\end{proof}

\begin{proof}[Proof of~\Cref{thmStairs}]
By the observation made at the beginning of~\Cref{sectionHooks}, we can assume $ d = 0 $. Setting $ m = \lambda_1+1-n $ in~\Cref{lemmaStairs}, we obtain
\[
	A_\lambda(n) = \sum_{i=1}^n\binom{\lambda_1+i-n}{\lambda_1+1-n}A_V(2n+1, i).
\]
Changing the index of summation from $ i $ to $ j = n-i+1 $ and using $ A_V(2n+1, i) = A_V(2n+1, n-i+1) $, we get
\[
	A_\lambda(n) = \sum_{j=1}^n \binom{\lambda_1+1-j}{\lambda_1+1-n}A_V(2n+1, j).
\]
Finally, it follows from~\Cref{thmRefinedVSASM} that
\[
	A_\lambda(n) = \frac{A_V(2n-1)}{\binom{4n-2}{2n-1}}\sum_{j=1}^n \binom{\lambda_1+1-j}{\lambda_1+1-n}\binom{2n+j-2}{2n-1}\binom{4n-j-1}{2n-1}.\qedhere
\]
\end{proof}

\bibliographystyle{abbrv}
\bibliography{StuNTSpapers}
\end{document}